\documentclass[dvips,coll,noinfoline]{article} %authoryear, ,

\RequirePackage[lnms]{imsart}
\RequirePackage{imsartprelims}
\RequirePackage[colorlinks]{hyperref}
\pubyear{2008}
\volume{3}
\volumetitle{Pushing the Limits of Contemporary Statistics:
Contributions in Honor of\\ Jayanta K. Ghosh}
\firstpage{1}
\lastpage{8}

\startlocaldefs
\setattribute{journal}{name}{Collections}
\setattribute{printed}{text}{Printed in Lithuania}
\endlocaldefs

\begin{document}

\begin{titlepage}
\editor{Bertrand Clarke and Subhashis Ghosal, Editors}
\end{titlepage}

\begin{copyrightpage}

\LCCN{2008924408}
%\ISBN[13]{???}
%\ISBN[10]{???}
\ISBN{978-0-940600-75-1}
\ISSN{1939-4039}
\serieseditor{Anthony Davison}
\treasurer{Rong Chen}
\executivedirector{Elyse Gustafson}

\end{copyrightpage}
%
%\begin{files}
%\filelist{files.txt}
%\end{files}
%
%\end{document}

\makeatletter
\gdef\doi@base{http://arXiv.org/abs/}
\makeatother

\begin{contents}[doi]
\contentsline{begintocitem}{}{}
\contentsline{jobname}{imscoll1pr}{}
\contentsline{doi}{0806.4445}{}
\contentsline{title}{Preface}{v}
\contentsline{author}{Bertrand Clarke and Subhashis Ghosal}{v}
\contentsline{endtocitem}{}{}

\contentsline{begintocitem}{}{}
\contentsline{jobname}{lnms55PR}{}
\contentsline{doi}{0806.4445}{}
\contentsline{title}{Contributors}{vii}
\contentsline{author}{}{vii}
\contentsline{endtocitem}{}{}

\contentsline{begintocitem}{}{}
\contentsline{jobname}{imscoll301}{}
\contentsline{doi}{0805.3066}{}
\contentsline{title}{J. K. Ghosh's contribution to statistics: A brief outline}{1}
\contentsline{author}{Bertrand Clarke and Subhashis Ghosal}{1}
\contentsline{endtocitem}{}{}

\contentsline{begintocitem}{}{}
\contentsline{jobname}{imscoll302}{}
\contentsline{doi}{0805.3064}{}
\contentsline{title}{Objective Bayesian analysis under sequential experimentation}{19}
\contentsline{author}{Dongchu Sun and James O. Berger}{19}
\contentsline{endtocitem}{}{}

\contentsline{begintocitem}{}{}
\contentsline{jobname}{imscoll303}{}
\contentsline{doi}{0805.3070}{}
\contentsline{title}{Sequential tests and estimates after overrunning based on $p$-value combination}{33}
\contentsline{author}{W. J. Hall and Keyue Ding}{33}
\contentsline{endtocitem}{}{}

\contentsline{begintocitem}{}{}
\contentsline{jobname}{imscoll304}{}
\contentsline{doi}{0805.3073}{}
\contentsline{title}{On predictive probability matching priors}{46}
\contentsline{author}{Trevor J. Sweeting}{46}
\contentsline{endtocitem}{}{}

\contentsline{begintocitem}{}{}
\contentsline{jobname}{imscoll305}{}
\contentsline{doi}{0805.3203}{}
\contentsline{title}{Data-dependent probability matching priors for empirical and related likelihoods}{60}
\contentsline{author}{Rahul Mukerjee}{60}
\contentsline{endtocitem}{}{}

\contentsline{begintocitem}{}{}
\contentsline{jobname}{imscoll306}{}
\contentsline{doi}{0805.3204}{}
\contentsline{title}{Probability matching priors for some parameters of the bivariate normal distribution}{71}
\contentsline{author}{Malay Ghosh, Upasana Santra and Dalho Kim}{71}
\contentsline{endtocitem}{}{}

\contentsline{begintocitem}{}{}
\contentsline{jobname}{imscoll307}{}
\contentsline{doi}{0805.3205}{}
\contentsline{title}{Fuzzy set representation of a prior distribution}{82}
\contentsline{author}{Glen Meeden}{82}
\contentsline{endtocitem}{}{}

\contentsline{begintocitem}{}{}
\contentsline{jobname}{imscoll308}{}
\contentsline{doi}{0805.3209}{}
\contentsline{title}{Fuzzy sets in nonparametric Bayes regression}{89}
\contentsline{author}{Jean-Fran\c {c}ois Angers and Mohan Delampady}{89}
\contentsline{endtocitem}{}{}

\contentsline{begintocitem}{}{}
\contentsline{jobname}{imscoll309}{}
\contentsline{doi}{0805.3220}{}
\contentsline{title}{Objective Bayes testing of Poisson versus inflated Poisson models}{105}
\contentsline{author}{M. J. Bayarri, James O. Berger and Gauri S. Datta}{105}
\contentsline{endtocitem}{}{}

\contentsline{begintocitem}{}{}
\contentsline{jobname}{imscoll310}{}
\contentsline{doi}{0805.3224}{}
\contentsline{title}{Consistent selection via the Lasso for high dimensional approximating regression models}{122}
\contentsline{author}{Florentina Bunea}{122}
\contentsline{endtocitem}{}{}

\contentsline{begintocitem}{}{}
\contentsline{jobname}{imscoll311}{}
\contentsline{doi}{0805.3238}{}
\contentsline{title}{Asymptotic optimality of a cross-validatory predictive approach to linear model selection}{138}
\contentsline{author}{Arijit Chakrabarti and Tapas Samanta}{138}
\contentsline{endtocitem}{}{}

\contentsline{begintocitem}{}{}
\contentsline{jobname}{imscoll312}{}
\contentsline{doi}{0805.3244}{}
\contentsline{title}{Risk and resampling under model uncertainty}{155}
\contentsline{author}{Snigdhansu Chatterjee and Nitai D. Mukhopadhyay}{155}
\contentsline{endtocitem}{}{}

\contentsline{begintocitem}{}{}
\contentsline{jobname}{imscoll313}{}
\contentsline{doi}{0805.3248}{}
\contentsline{title}{Remarks on consistency of posterior distributions}{170}
\contentsline{author}{Taeryon Choi and R. V. Ramamoorthi}{170}
\contentsline{endtocitem}{}{}

\contentsline{begintocitem}{}{}
\contentsline{jobname}{imscoll314}{}
\contentsline{doi}{0708.4294}{}
\contentsline{title}{Large sample asymptotics for the two-parameter Poisson--Dirichlet process}{187}
\contentsline{author}{Lancelot F. James}{187}
\contentsline{endtocitem}{}{}

\contentsline{begintocitem}{}{}
\contentsline{jobname}{imscoll315}{}
\contentsline{doi}{0805.3252}{}
\contentsline{title}{Reproducing kernel Hilbert spaces of Gaussian priors}{200}
\contentsline{author}{A. W. van der Vaart and J. H. van Zanten}{200}
\contentsline{endtocitem}{}{}

\contentsline{begintocitem}{}{}
\contentsline{jobname}{imscoll316}{}
\contentsline{doi}{0805.3264}{}
\contentsline{title}{A Bayesian semi-parametric model for small area estimation}{223}
\contentsline{author}{Donald Malec and Peter M\"uller}{223}
\contentsline{endtocitem}{}{}

\contentsline{begintocitem}{}{}
\contentsline{jobname}{imscoll317}{}
\contentsline{doi}{0805.3269}{}
\contentsline{title}{A hierarchical Bayesian approach for estimating the origin of a mixed population}{237}
\contentsline{author}{Feng Guo, Dipak K. Dey and Kent E. Holsinger}{237}
\contentsline{endtocitem}{}{}

\contentsline{begintocitem}{}{}
\contentsline{jobname}{imscoll318}{}
\contentsline{doi}{0805.3273}{}
\contentsline{title}{Kendall's tau in high-dimensional genomic parsimony}{251}
\contentsline{author}{Pranab K. Sen}{251}
\contentsline{endtocitem}{}{}

\contentsline{begintocitem}{}{}
\contentsline{jobname}{imscoll319}{}
\contentsline{doi}{0805.3279}{}
\contentsline{title}{Orthogonalized smoothing for rescaled spike and slab models}{267}
\contentsline{author}{Hemant Ishwaran and Ariadni Papana}{267}
\contentsline{endtocitem}{}{}

\contentsline{begintocitem}{}{}
\contentsline{jobname}{imscoll320}{}
\contentsline{doi}{0805.3282}{}
\contentsline{title}{Nonparametric statistics on manifolds with applications to shape spaces}{282}
\contentsline{author}{Abhishek Bhattacharya and Rabi Bhattacharya}{282}
\contentsline{endtocitem}{}{}

\contentsline{begintocitem}{}{}
\contentsline{jobname}{imscoll321}{}
\contentsline{doi}{0805.3286}{}
\contentsline{title}{An ensemble approach to improved prediction from multitype data}{302}
\contentsline{author}{Jennifer Clarke and David Seo}{302}
\contentsline{endtocitem}{}{}

\contentsline{begintocitem}{}{}
\contentsline{jobname}{imscoll322}{}
\contentsline{doi}{0805.3287}{}
\contentsline{title}{Sharp failure rates for the bootstrap particle filter in high dimensions}{318}
\contentsline{author}{Peter Bickel, Bo Li and Thomas Bengtsson}{318}
\contentsline{endtocitem}{}{}

\end{contents}

\begin{preface}

  \begin{frontmatter}

    \title{Preface}

  \end{frontmatter}

  \thispagestyle{plain}
Jayanta Kumar Ghosh is one of the most extraordinary professors in the field of Statistics.
His research in numerous areas, especially asymptotics, has been groundbreaking, influential
throughout the world, and widely recognized through awards and other honors. His
leadership in Statistics as Director of the Indian Statistical Institute and President
of the International Statistical Institute, among other eminent positions, has been
likewise outstanding. In recognition of Jayanta's enormous impact, this volume is
an effort to honor him by drawing together contributions to
the main areas in which he has worked and continues to work.  The papers
naturally fall into five categories.

First, sequential estimation was Jayanta's
starting point.  Thus, beginning with that topic, there are two papers, one
classical by Hall and Ding leading to a variant on p-values, and one Bayesian
by Berger and Sun extending reference priors to stopping time problems.

Second, there are five papers in the general area of prior specification. Much
of Jayanta's earlier work involved group families as does Sweeting's paper
here
for instance.   There are also two papers dwelling on the link between fuzzy
sets and priors, by Meeden and by Delampady and Angers. Equally daring is the work
by Mukerjee with data dependent priors and the pleasing confluence of several
prior selection criteria found by Ghosh, Santra and Kim.  Jayanta himself studied a
variety of prior selection criteria including probability matching
priors and reference priors.

Third, between his work on parametric Bayes and nonparametrics, Jayanta took an
interest in model selection.  Accordingly, three papers on model selection come
next.  Bunea's work on consistency echoes Jayanta's work on consistency of the
BIC.  Chatterjee and Mukhopadhyay's work on data adaptive model averaging continues the
direction they started under Jayanta's guidance.  Chakrabarti and Samanta's
work on the asymptotic optimality of
predictive cross validation contrasts nicely with standard Bayes model
selection, via the BIC for instance.

Fourth, there are five papers generally on Bayesian nonparametrics.  Some are applied
as in Malec and Mueller's work on semi-parametrics in small area estimation or Guo, Dey and Holsinger's
work carefully using prior selection for modeling purposes.  And some are more
theoretical: Choi and Ramamoorthi provide a review, with some new results, on
posterior consistency while James focuses on a class of priors and van
der Vaart and van Zanten focus on the role of reproducing
kernel Hilbert spaces in Bayesian nonparametrics with Gaussian process
priors.

Finally, Jayanta has most recently turned his attention to high dimensional
problems.  On this topic, there are five papers from a variety of standpoints.
For instance, it is
possible to make unexpected use of the information in the large dimensions
themselves as in Sen's work with Kendall's tau.  Others focus on the parametric
parts of a nonparametric model as in Ishwaran and Papana, or in Bhattacharya and
Bhattachcarya.  A third tack in Clarke and Seo is the focus on selecting the dimensions
for use in emerging model classes. Finally, the work of Bickel, Li and Bengtsson
establishes a general convergence result for computing conditional distributions.

As can be seen, some papers fit comfortably into more than one section and
some only fit into a section if it is interpreted broadly.  Even so, we would like to think
that the papers have achieved a nice tradeoff between clustering rather nicely
around the topic of each category and maintaining a reasonable diversity in
line with Jayanta's work.

Despite his manifold research interests, asymptotics and their applications have
been the main recurring theme of Jayanta's research since he published
his first paper in sequential statistics in 1960 (at the age of 23).  So, as a
generality, asymptotics undergirds most of the material in this volume
honoring him.

Fortunately, asymptotic thinking pervades statistical inference, even in the
most applied contexts, so, this is hardly a limitation. On the other hand,
asymptotics has a way of being impenetrably abstract.  However, all the papers
here, are, in Woodroofe's memorable phrase, written at a level that would
be `accessible to a determined graduate student'.  We encourage readers to have
a look at least at the introductions of papers outside their research area,
just for pure love of the field and the joy of intellectual stimulation.  We
suspect that once someone has read the introduction, he or she will be
ineluctably led to finish reading the whole paper.

As editors, we have been delighted at the depth and quality of work our
contributors submitted. They all make foundational points in the spirit of Jayanta.
We believe each paper will be of interest to researchers, theoretical
and applied, who confront problems that are difficult enough that conventional
solutions are inadequate and closed form solutions are intractable in the
several areas covered here. We are deeply grateful to all contributors for
offering their finest work to this volume.

Of course, no volume such as this could have been possible without the free and
anonymous labor of referees:  You folks know who you are, but for the sake
of confidentiality we cannot name you.  We especially thank those who provided extremely
prompt reports when we badly needed them. If any of you meet one of us
at a conference, we owe you a drink.  Probably two -- you helped us immeasurably.

In terms of actually producing this volume, Jennifer Clarke provided invaluable
support.  She helped us repeatedly with compiling complete versions of the volume.
In particular, the final, detailed copy-editing was largely her work; the balance of
her account in the Bank of Karma is astronomical.  We can't thank her enough.

Along the way there were many other people who gave us their time and expertise.
Dipak Dey helped guide us as we prepared initial proposal.  Rick Vitale, the former editor for
the LNMS series, also did a first rate job in explaining the details of how we had go about this kind of project,
if we wanted it to be successful.  Rick then gave the initial approval -- thanks, Rick!
%Sujit Ghosh also helped us on several occasions when we were in a jam; we are grateful
%to him for his willingness to find time for us on short notice.
Anthony Davison, the current editor for the LNMS series, worked closely with us to
get the volume in final form and then gave it final approval.  We appreciate
the burden that you carried for us, Anthony.

Finally, we are grateful that the IMS put its scarce resources into supporting
this volume.  The IMS has a tradition of honoring its most illustrious
members and, in our view, Jayanta is most assuredly among them.

It was a pleasure to put this tribute together and we hope that in some small
way we have served the interests of the research world.

\bigskip
{\flushright

\begin{tabular}{l@{}}
Bertrand Clarke\\
Subhashis Ghosal\vspace*{10pt}\\
November 2007
%Melbourne, Australia
\end{tabular}

}

\end{preface}

\newpage
\begin{contributors}
\begin{theindex}

  \item Angers, J.-F., \textit {Universit\'e de Montr\'eal, Montr\'eal, QC, Canada}

  \indexspace

  \item Bayarri, M. J., \textit {University of Valencia, Valencia, Spain}
  \item Bengtsson, T., \textit {Bell Labs, Murray Hill, NJ, USA}
  \item Berger, J. O., \textit {Duke University, Durham, NC, USA}
  \item Bhattacharya, A., \textit {University of Arizona, Tucson, AZ, USA}
  \item Bhattacharya, R., \textit {University of Arizona, Tucson, AZ, USA}
  \item Bickel, P., \textit {University of California, Berkeley, CA, USA}
  \item Bunea, F., \textit {Florida State University, Tallahassee, FL, USA}

  \indexspace

  \item Chakrabarti, A., \textit {Indian Statistical Institute, Kolkata, India}
  \item Chatterjee, S., \textit {University of Minnesota, Minneapolis, MN, USA}
  \item Choi, T., \textit {Inha University, Incheon, Korea}
  \item Clarke, B., \textit {University of British Columbia, Vancouver, BC, Canada}
  \item Clarke, J., \textit {University of Miami School of Medicine, Miami, FL, USA}

  \indexspace

  \item Datta, G. S., \textit {University of Georgia, Athens, GA, USA}
  \item Delampady, M., \textit {Indian Statistical Institute, Bangalore, India}
  \item Dey, D. K., \textit {University of Connecticut, Storrs, CT, USA}
  \item Ding, K., \textit {Queen's University, Kingston, ON, Canada}

  \indexspace

  \item Ghosal, S., \textit {North Carolina State University, Raleigh, NC, USA}
  \item Ghosh, M., \textit {University of Florida, Gainesville, FL, USA}
  \item Guo, F., \textit {Virginia Polytechnic Institute and State University, Blacksburg, VA, USA}

  \indexspace

  \item Hall, W. J., \textit {University of Rochester, Rochester, NY, USA}
  \item Holsinger, K. E., \textit {University of Connecticut, Storrs, CT, USA}

  \indexspace

  \item Ishwaran, H., \textit {Cleveland Clinic, Cleveland, OH, USA}

  \indexspace

  \item James, L. F., \textit {Hong Kong University of Science and Technology, Kowloon, Hong Kong}

  \indexspace

  \item Kim, D., \textit {Kyungpook National University, Taegu, Korea}

  \indexspace

  \item Li, B., \textit {Tsinghua University, Beijing, China}

  \indexspace

  \item M\"uller, P., \textit {M.D. Anderson Cancer Center, Houston, TX, USA}
  \item Malec, D., \textit {U.S. Census Bureau, Washington, DC, USA}
  \item Meeden, G., \textit {University of Minnesota, Minneapolis, MN, USA}
  \item Mukerjee, R., \textit {Indian Institute of Management Calcutta, Kolkata, India}
  \item Mukhopadhyay, N., \textit {Virginia Commonwealth University, Richmond, VA, USA}

  \indexspace

  \item Papana, A., \textit {Case Western Reserve University, Cleveland, OH, USA}

  \indexspace

  \item Ramamoorthi, R. V., \textit {Michigan State University, East Lansing, MI, USA}

  \indexspace

  \item Samanta, T., \textit {Indian Statistical Institute, Kolkata, India}
  \item Santra, U., \textit {University of Florida, Gainesville, FL, USA}
  \item Sen, P. K., \textit {University of North Carolina, Chapel Hill, NC, USA}
  \item Seo, D., \textit {University of Miami School of Medicine, Miami, FL, USA}
  \item Sun, D., \textit {University of Missouri, Columbia, MO, USA}
  \item Sweeting, T. J., \textit {University College London, London, UK}

  \indexspace

  \item van der Vaart, A. W., \textit {Vrije Universiteit, Amsterdam, The Netherlands}
  \item van Zanten, J. H., \textit {Vrije Universiteit, Amsterdam, The Netherlands}

\end{theindex}

\end{contributors}

%\newpage
%
%%\thispagestyle{empty}
%  \thispagestyle{plain}
%\begin{figure}
%  \includegraphics{imscoll1pic}
%\end{figure}

\end{document}